\documentclass[11pt,onecolumn]{article}
\usepackage{amscd}
\usepackage{times}
\usepackage{amsmath}
\usepackage{amssymb}
\usepackage{xspace}
\usepackage{theorem}
\usepackage{graphicx}
\usepackage{ifpdf}
\usepackage{url,hyperref}
\usepackage{latexsym}
\usepackage{euscript}
\usepackage{xspace}
\usepackage{color}
\usepackage{makeidx}
\usepackage{picins,wrapfig}

\long\def\remove#1{}

\newtheorem{theorem}{Theorem}[section] 

\newtheorem{definition}[theorem]{Definition}
\newtheorem{proposition}[theorem]{Proposition}
\newenvironment{proof}{{\em Proof:}}{\hfill{\hfill\rule{2mm}{2mm}}}

\def\marrow{{\marginpar[\hfill$\longrightarrow$]{$\longleftarrow$}}}

\newcommand{\Dan}[1] {{\sc Dan says: }{\marrow\sf [#1]}}

\newcommand {\mm}[1] {\ifmmode{#1}\else{\mbox{\(#1\)}}\fi}

\newcommand{\Z}                        {\mathrm {\mathbb{Z}}}

\newcommand{\XX}		{{\sf X}}
\newcommand{\YY}		{{\sf Y}}
\newcommand{\TT}		{{\sf T}}
\newcommand{\LL}		{{\sf L}}

\newcommand{\lx}[1]	{{\XX}_{#1}}
\newcommand{\subx}[1]   {{\XX}_{\downarrow#1}}
\newcommand{\rx}[2]	{{\XX}_{#1,#2}}

\newcommand{\cancel}[1]

\begin{document}

\title{Defining and Computing Topological Persistence for $1$-cocycles}

\author{
Dan Burghelea\thanks{
Department of Mathematics,
The Ohio State University, Columbus, OH 43210,USA.
Email: {\tt burghele@math.ohio-state.edu}}
\quad\quad
Tamal K. Dey\thanks{
Department of Computer Science and Engineering,
The Ohio State University, Columbus, OH 43210, USA.
Email: {\tt tamaldey@cse.ohio-state.edu}}
\quad\quad
Du Dong\thanks{
Department of Mathematics,
The Ohio State University, Columbus, OH 43210, USA.
Email: {\tt dudong@math.ohio-state.edu}}
}

\date{}
\maketitle

\begin{abstract}
The concept of topological persistence, introduced recently
in computational topology, finds applications in studying 
a map in relation to the topology of its domain. Since its introduction,
it has been extended and generalized in various directions.
However, no attempt has been made so far to extend the
concept of topological persistence to a generalization
of `maps' such as {\em cocycles} which are discrete counterparts
of closed differential degree one-forms, a well known concept in differential geometry.
We define a notion of topological persistence for $1$-cocycles
in this paper and show how to compute its relevant numbers.
It turns out that, instead of the standard persistence,
one of its variants which we call {\em level persistence} can be
leveraged for this purpose. It is worth mentioning that
$1$-cocycles appear in practice such as in data ranking or
in discrete vector fields.
\end{abstract}
\thispagestyle{empty}
\setcounter{page}{0}
\newpage

\section{Introduction}
The {\em persistent homology}, introduced recently~\cite{ELZ02} in 
computational topology, studies the longevity of topological
features as one sweeps a topological space with growing sublevel
sets of a real-valued function. Since its introduction, the concept
has been generalized in many directions, including its study
for homology groups under different coefficient rings~\cite{ZC05},
its stability under function perturbations~\cite{CEH07}, its
multidimensional version~\cite{CZ09}, and 
other variants~\cite{CSD09,CEH09,DW07,FL10}. 

In a discrete setting, which is
prevalent in practice, the topological space in question
is usually taken as the one defined by a simplicial complex $\XX$
whereas the function in question is taken as a linear map
$f\colon \XX\to \mathbb R$ whose restriction to each
simplex is linear. Such a pair $(\XX,f)$ provides what is called
a $0$-cochain in algebraic topology. A $0$-cochain can be specified 
by providing only a $n$-dimensional vector of $f$-values 
on $n$ vertices of $\XX$.
Such data appear frequently in applications involving shape and data analysis. 
Similarly, a more general input may specify real values
on edges of the simplicial complex by presenting
a skew symmetric $n\times n$ matrix 
for a simplicial complex $\XX$ with $n$ vertices.
Such a data, in reverence to the term
$0$-cochain, is called a $1$-cochain.
An interesting situation, which occurs in practice, is 
when the $1$-cochain is a $1$-cocycle (definition in section 2).
For example, in data ranking, pairwise orderings which are locally 
consistent represent a $1$-cocycle~\cite{JLY,YY}. 
A sampling of a curl-free vector field on a Riemannian manifold 
provides interesting examples of $1$-cocycles. 

It is natural to study  the qualitative complexity  of 
$0$-cochains 
and $1$-cocycles in analogy with
Morse theory and Novikov theory~\cite{Novikov} respectively.
In particular, one would like to relate 
the dynamics  of the $0$-cochain and $1$-cocycle  
with the topology of $\XX$. 
For $0$-cochains, this is precisely achieved by
studying persistent topology. Although the concept
has been generalized in many directions as we mentioned,
it has not yet been extended to $1$-cocycles.
It is appropriate to mention that
persistence and  circle valued maps are also in the attention of \cite{SJ09} 
where the existence of a circle valued map 
is established using standard persistence.

 

In this paper, we provide a definition of persistence for almost integral
$1$-cocycles\footnote{the general case can be 
reduced to almost integral case} defined on a simplicial 
complex and an algorithm
to compute relevant numbers. First, we show 
that a $1$-cocycle when almost integral 
\footnote {values on integral $1$-cycles are integer multiples 
of a fixed real number} 
can be interpreted as a
$\mathbb S^1$-valued (circle valued) map. 

The standard persistence~\cite{ELZ02,ZC05} can not be 
applied to circle valued map 
since sublevel sets are not defined. 
We observe that {\em level persistence} originally introduced
as {\em interval persistence} in~\cite{DW07} can be leveraged
to define a notion of persistence for circle  valued continuous maps. 
This persistence uses a lifting of the circle valued map
which becomes a $\mathbb R-$valued  map on a lifted space of $\XX$ which,
in topologists' terminology, is the infinite cyclic cover 
associated with the cohomology class of the almost integral cocycle. 
We show how one can calculate the relevant numbers for
level persistence for this
lifted map by adapting the algorithm for standard persistence.  
This adaptation is simplified by our observation that
the incidence structure of a {\em cell} complex involved
in level persistence computation can be derived easily 
from the given {\em simplicial} complex with a filtration.
It is appropriate to mention that the level persistence considered here
can be approached from the Zigzag persistence introduced
in~\cite {CSD09} based on quiver representations. 
This and the Morse Novikov theory 
for circle valued maps from the perspective of quiver 
representations will be topics of further research. 

{\bf Note}: After this work was done, we became aware that the 
``relevant numbers'' for the level persistence of an 
$\mathbb R$-valued map (as described later) can be derived from the invariants 
of the Zigzag persistence~\cite{CSD09}.
More recently we have noticed that they are actually equivalent.  
Appendix D of a more recent paper \cite{BD11} makes this precise.

\section {$1$-cocycle}
Let $\XX= ({\mathcal X}_0,{\mathcal X})$ be a simplicial complex where
${\mathcal X}_0$ is the set of
vertices, and  $\mathcal X$ is a set of finite subsets of ${\mathcal X}_0$ 
which satisfy the following properties :

\begin{enumerate}
\item ${\mathcal X}_0\subseteq \mathcal X$ and 

\item If $\sigma \in \mathcal X $ then  $\tau \subset \sigma $ $\Rightarrow$  
$\tau \in \mathcal X.$  
\end{enumerate}

The subsets in $\mathcal X$ of cardinality $k+1$  are called $k$-simplices. 
These subsets, considered
as ordered $k$-tuples up to an even permutation, represent oriented
$k$-simplices. We denote sets of oriented
$k$-simplices by ${\mathcal X}_k$.
We continue to denote by $\XX$ the canonical topological space 
(the geometric realization) 
associated with the combinatorial structure $\XX$.

Broadly
speaking, cochains are real valued functions that assign numbers
to oriented simplices.
We only consider $0$- and $1$-cochains which
are sufficient for this exposition. 
A $0$-{\em cochain}  is a function restricted
to vertices, $f: {\mathcal X}_0\to \mathbb R.$
A $0$-cochain can be  identified with a 
continuous map $f: \XX \to \mathbb R$  whose
restriction to each simplex is linear (a linear map). 
The  $0$-cochain  $f$ is {\em generic} if  $f: {\mathcal X}_0\to \mathbb R$ 
is injective.

%

The $1$-cochains are are maps  ${\bf f}: {\mathcal X}_1\to \mathbb R$ whose 
 domain is the set ${\mathcal X}_1$ of oriented edges of $\XX$. 
Recall that ${\mathcal X}_1$ can be regarded as  
the subset of 
${\mathcal X}_0 \times {\mathcal X}_0$ consisting of pairs 
$(x,y)$ which are vertices of a $1$-simplex.
The map $\bf f$ is a {\em $1$-cocycle}
if it satisfies :

\begin{enumerate}
\item $\bold f(x,y)= -\bold f(y,x)$ for any ordered pair 
$(x,y)\in {\mathcal X}_1$, 
and
\item if $ (x,y,z)\in {\mathcal X}_2$ then $\bold f(x,y) +   \bold
f(y,z) +   \bold f(z,x)=0;$
equivalently   $\bold f(x,y) +   \bold f(y,z) =   \bold f(x,z).$
\end{enumerate}

Clearly, a $0$-cochain $f$, or equivalently a linear map, 
provides a $1$-cocycle ${\delta f}$ defined by 
$${\delta f}(x,y)= f(y)-f(x).$$
Any $1$-cocycle $\bold f$ represents a cohomology class 
$<\bold f>\in H^1(\XX;\mathbb R)$ and any such cohomology class 
is represented by a $1$-cocycle. Two $1$-cocycles $\bold f_1$ 
and $\bold f_2$ represent the same cohomology class iff 
$\bold f_1-\bold f_2= {\delta f}$ for some $0$-cochain $f$.

An  {\it almost integral} $1$-cocycle is a pair $({\bf f},\alpha)$
where $\bf f$ is a $1$-cocycle whose values on integral $1$-cycles
are integer multiple of a fixed positive real
$\alpha$.\footnote{in the language of 
algebraic topology the cohomology class $<\bold f>\in H^1(\XX; \mathbb R)$ 
has degree of rationality $1$, or equivalently the image of the 
induced homomorphism 
$H_1(\XX; \mathbb Z)\to \mathbb R$ is $\alpha \mathbb Z\subset \mathbb R$.}

If $St(x)$ denotes the  star of the vertex 
$x\in {\mathcal X}_0$ (the star of any simplex is a
sub complex), a $1$-cocycle $\bold f$
defines a unique function 
$f_x: St(x)\to \mathbb R $ by the
formulae  $f_x (x)=0$ and
$f_x (y)= \bold f (x, y)$ for any vertex 
$y\neq x$ in $St(x).$ 
Clearly, 
$(f_x- f_y) (z)$ is constant  for any $z$ 
in a connected component of $St(x)\cap St(y)$.

Thus, a $1$-cocycle can be thought as a collection of 
linear maps $\{f_x: St(x)\to \mathbb R\}$  for each vertex $x$, such that the 
difference $f_x-f_y$ is
constant on  each connected component 
of $St(x)\cap St(y).$ Therefore, one can regard 
$\bold f\equiv \{f_x, x\in
{\mathcal X}_0\}$  as a {\em multivalued linear map}.
A 1-cocycle $\bold f$ is {\em  generic}  if all linear maps $ f_x$ 
are  generic, i.e., injective when restricted to vertices of $St(x)$.

\subsection{$1$-cocycles and circle valued maps}
\label{cocycle-circle-sec}
Consider a continuous {\em circle valued} map $f: \XX\to {\mathbb S}^1$.
Let $p: \mathbb R\to {\mathbb S}^1={\mathbb R}/{\alpha \mathbb Z}$ 
be the  map defined  by $p(t)= t $(mod $\alpha$), 
$\alpha$ a positive real number. For
any simplex $\sigma \in \XX$, the restriction  
$f|_{\sigma}$ admits liftings $\hat
f:\sigma \to \mathbb R$, i.e., $\hat f$ is a 
continuous map which satisfies $p\cdot
\hat f= f|_\sigma$.
The map $f: \XX\to \mathbb S^1$ is called  {\em linear} 
if, for any simplex $\sigma $, at least
one  of the liftings (and then any other) is  linear.
Any linear map $f: \XX\to \mathbb S^1= \mathbb R/\alpha \mathbb Z$  defines an
almost integral $1$-cocycle $({\bf f}, \alpha)$
by assigning to each edge $e$ with vertices $x,y$,
$$\bold f(x,y)=\hat f( y)- \hat f(x)$$ 
where $\hat f$ is a lift of $f$ restricted to $e$.
One can show the opposite, that is,
\begin{itemize}
\item an almost integral $1$-cocycle $(\bold f:{\mathcal X}_1\to \mathbb{R},\alpha)$ 
can be associated to a linear map 
$f: \XX\to \mathbb S^1= \mathbb R/\alpha \mathbb Z$.
\end{itemize}

\paragraph{ Covering associated with  a cohomology class 
$<{\bf f}>\in H^1(\XX; \mathbb R):$}

Regard $\XX$ as a topological space. Choose a base point  $x\in \XX$ 
and consider the space of continuous paths $\
\gamma:[0,1]\to \XX$ with $\gamma(0)=x,$ equipped with
compact open topology.
Make two continuous paths $\gamma_1$ and $\gamma _2$ equivalent iff 
the $\gamma_1 (1)= \gamma_2(1)$ and
the closed path $\gamma_1\star \gamma_2^{-1}$  
satisfies $<{\bf f}> ([\gamma_1\star \gamma_2^{-1}])=0.$ Here
$\star$ denotes the concatenation of the paths $\gamma_1$ and $\gamma^{-1}_2$
defined by $\gamma^{-1}_2(t)= \gamma_2(1-t),$  
and $[\gamma_1\star \gamma_2^{-1}]$ denotes the homology class of
$\gamma_1\star \gamma_2^{-1}.$ 
The quotient space $\tilde \XX,$ 
whose underlying set is the set of equivalence classes of paths,  
is equipped with
the canonical map $\pi:\tilde \XX\to \XX$ induced by assigning  
to $\gamma $ the point
$\gamma(1)\in \XX.$ The map $\pi$  is a local homeomorphism 
and $\tilde \XX$ is the
{\em total space} of a principal covering 
with group  $G$ where
$$G={\rm img}(<\bold f> : H_1 (\XX;\mathbb Z)\to \mathbb R).$$ 
When  $\bold f$ is almost integral, $G$ is isomorphic to 
$\mathbb Z.$
If $\XX$ is equipped with a triangulation, then 
$\tilde \XX$ gets  a triangulation 
whose simplices, when viewed as subsets of $\tilde \XX$ are
homeomorphic by  $\pi$ to simplices of $\XX ($ when viewed as subsets 
of $\XX$).

\paragraph{Construction of $\tilde f$ and $f$.}
We construct the circle valued map $f$ via its 
cyclic cover $\tilde f$ as follows.
\begin{enumerate}
\item[Step 1.] Consider $\pi:\tilde \XX\to  \XX$  the principal 
$\mathbb Z-$covering associated
with the cohomology class $<\bold f>$ defined by $\bold f.$  
This means that $\tilde \XX$ is a simplicial complex equipped 
with a free simplicial action $\mu: \mathbb Z\times \tilde \XX\to \tilde \XX$ 
whose quotient  space,
$\tilde \XX/ \mu,$ is the simplicial complex $\XX$.
The map $\pi $ identifies
to the map
$\tilde \XX\to \tilde \XX/ \mu = \XX$ and satisfies $\pi^\ast (<\bold f>)=0.$

Choose a vertex $x$ of $\XX$ and call it
a {\em base point}.
Notice that the vertices $\tilde{\mathcal X}_0$
of $\tilde \XX$ can be also described  
as equivalence classes of
sequences 
$\{x= x_0, x_1, \cdots x_{N-1}, x_N\}$ 
with $x_i$s being consecutive vertices of $\XX$ (i.e. $x_i,x_{i+1}$ are
vertices of an edge). Two such sequences,
$\{x=x_0,x_1,\cdots,x_{N-1},x_N\}$ 
and $\{x= y_0, y_1, \cdots y_{L-1}, y_L\}$  
are equivalent if $x_N=y_L$ and
the sequence
$\{x= z_0, \cdots , z_{N+L}= x\}$ with $z_i= x_i$  
if $i\leq N$ and  $z_{j+N}= y_{L-j}$
if $j\leq L,$ satisfies
$$\sum_{0\leq i\leq L+N-1} \bold f(z_i, z_{i+1})= 0.$$

\item[Step 2.] 
Define the map $\tilde f: \tilde{\mathcal X}_0 \to \mathbb R$  
by  $\tilde f(\tilde y):= \sum_{0\leq i\leq {L-1}} \bold f(y_i, y_{i+1})$
where $\tilde y\in \tilde{\mathcal X}_0$ is the vertex
corresponding to the equivalent class of $\{x=y_0,\cdots,y_L\}$.
The description of $\tilde \XX$ given above guarantees 
that $\tilde f$ is well defined. Extend $\tilde f $ 
to a linear map $\tilde f: \tilde \XX\to \mathbb R.$
Observe that if $\tilde y_1$ and $\tilde y_2$ satisfy 
$\pi(\tilde y_1)= \pi(\tilde y_2)$ then $\tilde f (\tilde y_1)- 
\tilde f (\tilde y_2) \in  \alpha \mathbb Z.$ In addition  if $\tilde e_1$ 
and $\tilde e_2$ are two edges of $\tilde \XX$ from $\tilde y_1$ to $\tilde
y_1'$ and $\tilde y_2$ to $\tilde y'_2$ respectively with 
$\pi (\tilde e_1)= \pi(\tilde e_2),$ 
then $\tilde f(\tilde y'_1)- \tilde f(\tilde y'_2)= \tilde f(\tilde y_1)- \tilde f(\tilde y_2).$ 
This implies that if $\tilde \sigma_1$ and
$\tilde\sigma_2$ are two simplices with $\pi(\tilde \sigma_1)= \pi(\tilde
\sigma_2)= \sigma$  and $\pi_i'$s  are the restrictions of
$\pi$ to $\tilde \sigma_i$  ($\pi$ are  bijections on their image), then 
 $\tilde f \cdot \pi^{-1}_1 -\tilde f \cdot \pi^{-1}_2:\sigma\to \mathbb R,$ 
is constant and this constant is an integer multiple of the
fixed real number $\alpha.$ 

\item[Step 3.]  
\parpic[1]{
\begin{minipage}{0.3\textwidth}
\[
\begin{CD}
\tilde \XX                           @>\tilde f>> \mathbb R\\
@V\pi VV                                                       @Vp VV\\
\tilde \XX/\mu = \XX @>f >>            \mathbb S^1
\end{CD}
\]
\end{minipage}
}
Observe that  the map 
$p\cdot \tilde f :\tilde \XX\to \mathbb R$  (with $p:\mathbb
R\to \mathbb S^1= \mathbb R/ \alpha \mathbb Z$)  factors through $\tilde \XX/\mu= \XX $ 
inducing a map from $\XX$ to $\mathbb S^1.$  This is our linear map $f$
whose associated $1$-cocycle  is $\bold f.$ 
The relation between 
$f$ and $\tilde f$ can be summarized by the commutative diagram on left.
\end{enumerate}

\section {Definitions of persistence}
This section defines persistence for almost
integral $1$-cocycles via circle valued maps. 
First, we expose the definitions for level
persistence in analogy to standard persistence~\cite{EH,ZC05},
and then extend them to circle valued maps. 

Let $f: \XX\to \mathbb R$ be a continuous map. 
For simplicity suppose $\XX$ is a nice compact
space and $f$ a nice map. This means that the 
homology $H_r(\cdot)$, $r\geq 0$,
of all levels and sub levels  are   finitely generated.
Define 
\begin{eqnarray*}
\subx{t} &:=& f^{-1}((-\infty,t]),\\
\lx{t} &:=& f^{-1}(t),~\mbox{ and  }\\
\rx{t_1}{t_2} &:=& f^{-1}([t_1, t_2]), t_2\geq t_1.
\end{eqnarray*}

\paragraph{Persistence of a $\mathbb R$-valued map.}
Standard {\em persistent homology} with coefficients 
in a {\em field} $\kappa$ ~\cite{ZC05} 
is the collection of vector spaces 
$B_r(t,t'):= {\rm img} (H_r(\subx{t})\to H_r(\subx{t'})),\ t\leq t'$.
It is understood that the homology $H_r$ is considered with coefficients in the field $\kappa.$  The corresponding Betti numbers are $\beta_r(t, t'):=\dim B_r(t,t').$ 
Equivalently,  one  can define persistent homology as the 
collection of vector spaces 
$K_r(t):= H_r(\subx{t})$ 
each of which is equipped with a filtration provided by subspaces:

$$
K_r(t,t'):= \ker(H_r(\subx{t})\to H_r(\subx{t'}))
\mbox{ with } K_r(t,t')\subseteq K_r(t, t'')\subseteq K_r(t)
\mbox{ for } t\leq t' \leq t''.
$$
The corresponding  Betti numbers are 
$\kappa_r(t):=\dim K_r(t)$ and 
$\kappa_r(t,t'):= \dim K_r(t,t').$
The two type of informations are equivalent. Indeed,
$B_r(t,t)= K_r(t)$ giving $\beta_r(t,t)= \kappa_r(t)$ and 
$B_r(t,t')$ is isomorphic to $ K_r(t)/ K_r(t,t')$ implying
$\beta_r(t,t')= \kappa_r(t)- \kappa_r(t,t').$

\begin{definition} Let $c \in H_r(\subx{t})\setminus 0.$ 
One says that 
\begin{itemize}
\item[(i)] $c$  is
born at $t',t'\leq t$, if 
$c\in {\rm img}\,(H_r(\subx{(t'+\epsilon)}) \to H_r(\subx{t}))$ but
$c\notin {\rm img}(H_r(\subx{(t'-\epsilon)}) \to H_r(\subx{t}))$ 
for  any positive $\epsilon$ 
with $t-t' >\epsilon > 0$,
\item[(ii)] $c$  dies at $t''$, $t'' \geq t,$
if its image is nonzero in ${\rm img}\,(H_r(\subx{t})\to 
H_r(\subx{(t''-\epsilon)}))$ 
but zero in  ${\rm img}\,(H_r(\subx{t})\to H_r(\subx{(t''+\epsilon)}))$ 
for any positive $\epsilon$  
with $t''-t >\epsilon > 0$. 
\end{itemize}
\end{definition}

For any $c\in H_r(\subx{t})$, let $t^+(c)$ and $t^-(c)$
denote the death and birth times of $c$ respectively where
$t^+(c) \geq t$, $t^-(c)\leq t.$
For $t'\leq t\leq t'',$
one may introduce the numbers 
\begin{itemize}
\item
$\mu_r(t', t, t'') :=$ {\em the maximal number of linearly independent 
elements} in  $H_r(\subx{t})$ so that  $t^+(x)= t''$ and $t^-(x)= t'.$ 
This also means that no linear combination of such elements   
is born before  $t'$, or dies before $t''$.
\end{itemize}

The numbers $\mu_r(t', t, t'')$ 
do not depend on $t,$  hence can be written as $\mu_r(t', t'')$.
The numbers  $\beta_r(t', t'')$ determine and are determined 
by the numbers $\mu_r(t', t''),$ cf \cite {EH}.
The set of numbers $\beta_r(t, t')$ or 
$\kappa_r(t), \kappa_r(t,t')$ or $\mu_r(t,t')$  
are the {\em relevant persistent numbers}.

\paragraph{Simultaneous persistence.}
Suppose $f_\pm: \XX_\pm \to [0,\infty)$ are two continuous maps as 
above with $f_+^{-1}(0)= f_-^{-1}(0)= A.$ It will be
useful to calculate the number
\begin{itemize}
\item
 $\omega_r(s, t) :=$ {\em the maximal number of linearly independent 
elements} in  $H_r(A)$  which die in $\XX_-$ at $s$ and in 
$\XX_+$ at $t$ (this  means that no linear combination of such elements   
dies before $s$ resp. $t$) is referred as the simultaneous persistence 
number of the pairs $\{f^\pm : \XX^\pm\to [0,\infty)\}.$
\end{itemize}
When $\XX$  is a cell complex and $f$  is a linear map restricted to each cell,
there are effective algorithms to compute the numbers $\mu_r(t', t'').$ 
These algorithms  can be adapted to compute the number $\omega_r(s,t)$ 
for a pair $f^\pm : \XX^\pm\to [0,\infty)$ 
with $\XX^\pm$ and $A$ being cell complexes and 
$f^\pm$ being linear on each cell
as indicted in section 4.

\subsection{Level Persistence}
\label{lpers-sec}
Level sets instead of sublevel sets define level persistence.
Let
\begin{eqnarray*}
L^+_r(t;\tau)&:=& \ker (H_r(\lx{t})\to H_r(\lx{t,t+\tau})) \mbox{ and }\\
L^-_r(t;\tau)&:=& \ker (H_r(\lx{t})\to H_r(\lx{t-\tau, t})).
\end{eqnarray*}

\begin{definition}
The {\em level persistent homology} with coefficients in 
a field  is the collection of vector spaces 
$L_r(t):= H_r(\lx{t})$ equipped with two filtrations:
\begin{eqnarray*}
L^+_r(t;\tau) \subseteq  L^+_r(t;\tau') \subseteq L_r(t), \tau \leq \tau'
\mbox{ and }\\
L^-_r(t;\tau) \subseteq  L^-_r(t;\tau') \subseteq L_r(t), \tau \leq \tau'.
\end{eqnarray*}
\end{definition}
Consequently we have the {\em relevant level persistence
numbers}: 
\begin{eqnarray*}
l_r(t):= \dim L_r(t),~ l_r^\pm (t;\tau) :=\dim L^\pm _r(t;\tau),
\mbox{ and }\\
e_r(t; \tau',\tau''):= \dim (L^-_r(t;\tau')\cap L^+_r(t;\tau'')).
\end{eqnarray*}

\begin{definition} Let $c \in H_r(\lx{t})$. One says that 
\begin{itemize}
\item[(i)] $c$ dies 
downward at $t'  ,t' \leq t$, if its image is nonzero in 
${\rm img}\,(H_r(\lx{t})\to H_r(\rx{t'+\epsilon}{t}))$ 
but is zero in ${\rm img}\,(H_r(\lx{t})\to H_r(\rx{t'}{t}))$ for  
any $\epsilon$ with $t-t'>\epsilon>0$, 
\item[(ii)]
$c$ dies upward  at $t'',t''\geq t$, if its image 
is nonzero in ${\rm img}\,(H_r(\lx{t})\to H_r(\rx{t}{t''-\epsilon}))$ but 
is zero in ${\rm img}\,(H_r(\lx{t})\to H_r(\rx{t}{t''}))$ for 
any $\epsilon$ with $t''-t >\epsilon > 0.$
\end{itemize}
\end{definition}

For $c\in H_r(\lx{t})$, let $\tau^+(c)$ and  $\tau^-(c)$
be the upward and 
downward {\em life} time of $c$ respectively, 
if $c$ dies upward and downward at 
$t+\tau^+(c)$ and $t- \tau^-(c)$ respectively. 
For $\tau^-, \tau^+\geq 0$, and $t$, we introduce the numbers :
\begin {itemize}
\item 
$\nu^+_r( t; \tau^+) :=$ the maximal number of linearly independent 
elements in 
$H_r(\lx{t})$ so that  $\tau^+(c)= \tau^+,$ 
\item 
$\nu^-_r( t; \tau^-) :=$ the maximal number of linearly independent elements in 
$H_r(\lx{t})$ so that  $\tau^-(c)=\tau^-.$

\end{itemize}
For tame maps, the numbers  $\l^\pm_r(t; \tau)$  
determine and are determined by the numbers $\nu^\pm_r (t;\tau)$.
The numbers $\nu^\pm_r(t;\tau)$
can be calculated from standard persistence
numbers $\mu_r (\cdot,\cdot) $ 
for another space $\YY$ equipped with a tame map $g:\YY\to \mathbb R$ 
derived from $(\XX, f:\XX\to \mathbb R).$  

Similarly, the numbers $e_r(t; \tau', \tau'')$ can be 
derived from  $\omega_r(s,t)$ calculated for a pair $( \YY^\pm, g^\pm)$
derived canonically from $(\XX, f)$.    
Precisely, 
\begin{enumerate}
\item $\nu^+_r(t; \tau^+)$ for map $f$ is same as $\mu _r(t, t+\tau^+)$
for map $g:\YY\to \mathbb{R}$ where $\YY=\rx{t}{\infty}$ 
and $g= f|_{\rx{t}{\infty}}$.

\item $\nu^-_r(t; \tau^-)$ for map $f$ is same as $\mu_r
(-t, - t+\tau^-)$ for map $g:\YY\to \mathbb R$ 
where $\YY= \rx{-\infty}{t}$ and $g= -f|_{\rx{-\infty} {t}}$.
\item $e_r(t; \tau^-, \tau^+)$ is the same as the  simultaneous persistence 
number $\omega_r(\tau^-,\tau^+)$ associated with $\YY^+= f^{-1}([t,\infty))$,
$\YY^-= f^{-1}((-\infty, t])$, $g_+= f|_{\YY_+}- t$,  and $g_-= -f |_{\YY_-} +t.$
\end{enumerate} 

\subsection{Persistence for circle valued maps}
The standard (sub level) persistence can not be extended to 
circle valued maps because the notion of
sub level is not well defined for circle valued maps. However, we can
extend the level persistence to such maps.

Let $f: \XX\to \mathbb S^1=\mathbb R/\alpha\mathbb Z$ be a 
continuous circle valued map.
Consider 
\begin{itemize}
\item $p: \mathbb R\to \mathbb S^1,$ defined by $p(t)= t (\rm{mod}~\alpha)$.

\item $\pi =\pi^f:\tilde \XX\to \XX,$ the $\mathbb Z-$ principal covering   
associated to $f$,  
i.e., the pull back of the principal covering $\mathbb R\to \mathbb S^1$ by the 
map $f$. Equivalently, it is the principal covering associated to 
the $1$-cocycle defined by $f$
as described in section~\ref{cocycle-circle-sec}  

\item $\tilde f: \tilde \XX\to \mathbb R$ is a map which 
satisfies $p \cdot \tilde
f= f\cdot\pi^f$. We make it unique by 
setting $\tilde f(\tilde x)=0$ for a chosen base point $\tilde x$ 
of $\tilde \XX.$
\end{itemize}
The map $\tilde f$ is the infinite cyclic 
covering of $f$ which would be same as the one constructed 
(section~\ref{cocycle-circle-sec}) from
the $1$-cocycle associated to $f$.
Observe  that $\pi(\tilde{\XX}_t)=  \lx{p(t)}= f^{-1}(p(t)).$

\begin{definition}
Define $L^\pm_r(\theta; \tau)$ for $f$ to 
be $L^\pm_r(t; \tau)$  for the
map $\tilde f$ where $p(t)= \theta$.
The definition involves the choice of $t$. It is easy to verify the
independence of this choice.
\end{definition}

Now, we are ready to define persistence for circle valued maps.
Let $\theta\in \mathbb S^1$ denote an angle and $\XX_{\theta}=f^{-1}(\theta)$.
The  {\em level persistent homology} with coefficients in 
a field  for $f$ 
is  the collection of vector spaces 
$L_r(\theta):= H_r(\lx{\theta})$ equipped with two filtrations:

\begin{itemize}
\item 
$L^+_r(\theta ; \tau) \subseteq  L^+_r(\theta;\tau') \subseteq L_r(\theta), 
0 \leq \tau \leq \tau'\leq \infty$,
\item
$L^-_r(\theta ;\tau) \subseteq  L^-_r(\theta;\tau') \subseteq L_r(\theta), 
0\leq \tau \leq \tau'\leq \infty$.
\end{itemize}

Consequently we have the relevant persistence numbers: 
\begin{eqnarray*}
l_r(\theta)&:=& \dim L_r(\theta),\\
l_r^\pm (\theta; \tau) &:=&\dim L^\pm _r(\theta; \tau), \mbox{ and }\\
e_r(\theta; \tau',\tau'')&:=& \dim (L^-_r(\theta; \tau')\cap 
L^+_r(\theta; \tau'').
\end{eqnarray*}

By the above definition, the numbers $l_r^\pm (\theta; \tau)$ 
and $e_r(\theta; \tau',\tau'')$ for $f$ are the numbers $l_r^\pm (t; \tau)$ 
and $e_r(t; \tau',\tau'')$ for $\tilde f.$ 
Note that the  numbers $l_r^\pm (\theta; \tau)$ and
$e_r(\theta; \tau' \tau'')$ can be computed from 
the level persistence numbers of $\tilde f$    
which in turn 
can be derived from persistence numbers  
for the maps $g$ and simultaneous persistence numbers for the 
pair $g^\pm$ as described in subsection~\ref{lpers-sec}.

\subsection{Persistence for almost integral 1-cocycles}
Given an almost integral 1-cocycle  $(\bold f,\alpha)$, 
we associate the circle ($\mathbb S^1= \mathbb R / \alpha \mathbb Z$) 
valued map 
$f: \XX\to \mathbb S^1$ as described in section~\ref{cocycle-circle-sec}
and define the persistence of $\bold f$ as the 
level persistence of the circle valued map $f$. 
It is possible to shortcut the involvement of the map $f$ 
and go directly to $\tilde f$ to define persistence of $\bf f$,
but this might  obscure the topological meaning of the definition. 



\subsection {Tame maps} 
\label{tame-sec}

For finite calculations, one cannot allow the level sets change
topology continuously. This is why we introduce the following
restrictions of tameness on the maps. Similar conditions
for standard persistence have been proposed before~\cite{CEH07}.
In most practical situations, the tameness
condition holds for the maps of interest. 

\begin{definition}.
A continuous map $f: \XX\to \YY,\  \YY=\mathbb R\  \rm {or}\   \mathbb S^1$ 
is called tame if there exists finitely many $\{t_1, t_2,\cdots t_N\}$ so that :
\begin{itemize}
\item[(i)] for any $t \ne t_1, t_2, \cdots t_N$ there exists $\epsilon>0$ 
so that $f: \lx{t-\epsilon, t+\epsilon} \to [t-\epsilon, t+\epsilon]$ 
and the second factor projection 
$\lx{t}\times [t-\epsilon, t+\epsilon] \to [t-\epsilon, t+\epsilon]$ 
are fiberwise homotopy equivalent,
\item[(ii)] for any $t$ (in particular for $t=t_i$) there exists $\epsilon >0$ so that the canonical inclusion $\lx{t} \subset \lx{t-\epsilon, t+\epsilon}$ is a homotopy equivalence.
\end{itemize}
\end{definition} 

The above definition can be considerably weakened 
by considering homology equivalence in place of homotopy
equivalence, but in view of the fact that all our examples 
satisfy the above definition, we proceed with it. 
Generic smooth maps from a closed smooth manifold to $\mathbb R$ or $\mathbb S^1$ are tame and so are piecewise linear maps from a simplicial complex to $\mathbb R$ 
or $\mathbb S^1.$ If $f:\XX\to \mathbb S^1$ is tame, 
then  $\tilde f: \tilde \XX \to \mathbb R$ 
restricted to any $\tilde{\XX}_{t_1, t_2}$ is tame.  

Observe that, for  a tame map $f :\XX \to \mathbb R$, 
the vector spaces $B_r(t,t')$  are completely determined by 
the vector spaces $B_r(i,j):= B_r(t_i,t_j)$ and therefore 
$\beta_r(t,t')$ and $\mu_r(t', t, t'')$ by the numbers 
$\beta_r(i,j):= \beta_r(t_i, t_j)$ and 
$\mu_r (i, j, k):= \mu_r(t_i, t_j, t_k)$.
The following relations between these numbers are well known~\cite{CEH07,EH}. 
\begin{eqnarray*}
\mu_r(i,k)= \mu_r(i, j, k)=  \beta_r (i, k-1)-\beta_r (i,k)-\beta_r(i-1, k-1) +
\beta_r(i-1, k),\\
\beta_r(i,j)= \sum _{j'>j,\  i'\leq i}\mu_r(i', j') \mbox{ with }
\dim H_r(\subx{t_i})= \beta_r(i,i).
\end{eqnarray*}

One may observe similar properties for level persistence. Let $s_{2i}= t_i$ 
and $s_{2i-1}$ be any number between $t_{i-1}$ and $t_i.$
Clearly, the vector spaces $L_r(t)$ and $L^\pm_r(t;\tau)$  are completely 
determined by the vector spaces:
 $L_r(i):= L_r(s_i)$, 
$L^+_r(i;k)= L^+_r(s_i; s_{i+k}- s_i),$
$L^-_r(i;k)= L^+_r(s_i; s_{i}- s_{i-k})$.
Therefore, the numbers 
$\l^\pm_r(t;\tau)$, $\nu_r^\pm(t;\tau)$, and $e_r (t;\tau^-,\tau^+)$ 
are determined by the numbers 

\begin{eqnarray*}
l^+_r(i; k):=& l^+_r(s_i; s_{i+k}-s_i)\\
l^-_r(i; k):= &l^-_r(s_i; s_i-s_{i-k})\\
\nu^+_r(i;k):=&\nu^+_r(s_i; s_{i+k}-s_i)\\
\nu^-_r(i;k):=&\nu^-_r(s_i; s_i-s_{i-k})\\
e_r(i; j,k):=&e_r(s_i; s_i-s_{i-j}, s_{i+k}-s_i) 
\end{eqnarray*}

Observe that 
\begin{eqnarray*}
&l^+_r(i;0)= l^-_r(i;0)=0\\
&l^+_r(2i;1)= l^-_r(2i;1)=0\\
&e_r(i;0,0)= e_r(i;0,1)= e_r(i;1,0)= e_r(2i;1,1)=0 
\end{eqnarray*}
and

$$l^+_r (i;j)= \sum_{0 \leq k\leq j} \nu ^+_r(i;k),~~
l^- _r(i;j)=\sum_{0\leq k \leq j} \nu ^-_r(i;k).$$

This shows that level persistence numbers can be computed from
$\nu^\pm_r(\cdot,\cdot)$ which in turn can be computed from
$\mu_r(\cdot,\cdot)$ (section~\ref{lpers-sec}).
Similarly, the level persistence numbers $e_r(i; j,k)$ 
can be calculated by the simultaneous persistence numbers
$\omega_r(s_i- s_{i-j}, s_{i+k}-s_i) $
 by the formula $$e_r(i; j,k)= \sum_{0\leq j'\leq j;
 0\leq k'\leq k} \omega_r(s_i-s_{i-j'}, s_{i+k'}- s_i)$$

\section {Algorithm}
The algorithm to compute persistence
for an almost integral $1$-cocycle $\bf f$ follows the
logical sequence that defines its relevant persistence
numbers in the previous section. It considers 
the associated circle valued map $f$, and
computes the level persistence for the cyclic 
covering $\tilde f$. We assume $\Z_2$-homology
and describe how one can adapt a matrix version~\cite{CEM06,EH} 
of the standard persistence algorithm~\cite{ELZ02} to compute 
level persistence for a $\mathbb R$-valued map,
then how to extend it to circle valued maps.
In this development, we encounter subspaces
of simplicial complexes which are cell complexes.
We begin with its definition.

\begin{definition}
A geometric cell complex  $\XX$ is a union of  a 
collection $\mathcal X$ of 
non degenerate convex cells with disjoint interiors 
which satisfy the property that, for any cell $\sigma\in \mathcal X$, 
all its faces belong to $\mathcal X$.
\end {definition}
We consider finite cell complexes $\XX$ equipped with a 
filtration $ \mathcal F:\equiv \XX_0\subseteq 
\cdots  \XX_k \subseteq \XX_{k+1}\cdots \subseteq \XX_m=\XX$ 
with $\XX_i$ being sub complexes of $\XX$. 

\begin{definition}
A total order on $\mathcal X=\{\sigma_1,\ldots,\sigma_n\}$
is called topologically consistent if the condition A below is satisfied 
and filtration compatible if the condition B below is satisfied.
\begin {itemize}
\item
{\bf Condition A.} 
$\sigma_i$ is a face of $\sigma_j$ implies $i< j.$ 
\item
{\bf Condition B.}
$\sigma_i \in \XX_k$ and $\sigma_j \in \XX_{k'}\setminus \XX_k$ with $k<k'$ 
implies $i<j.$
\end{itemize}
\end{definition}
Given a filtration $\mathcal F$, one can canonically  modify any total order 
which satisfies Condition A into one which satisfies 
both conditions A and B. 
This is done in the following way:\\

Find the first simplex $\sigma$ in the order which violates Condition B. 
Let $\tau$ be the simplex immediately preceding $\sigma$ in the total order.
The violation of condition B by $\sigma$ implies that
$\tau\in \XX_j\setminus \XX_{j-1}$ and $\sigma\in \XX_i\setminus \XX_{i-1}$ where $i<j$. Permute $\sigma$ with $\tau$.
Observe that condition A continues to hold; indeed the only possible 
violation of condition A is if $\tau$ is a face of $\sigma$ 
which would imply that $j\leq i$.  If $\sigma$ still violates Condition B
in the new position, then  continue moving $\sigma$ to left
until it does  not violate Condition B anymore. All other simplices which 
were initially preceding $\sigma$ do satisfy condition B.
So, we have one less simplex which violates Condition B.

\paragraph{Persistence algorithm~\cite{CEM06,ELZ02}.}
The input to the algorithm is a matrix $M$ that represents
incidence structure of a cell complex $\XX$ equipped with a filtration. 
The algorithm derives the numbers $\mu_r(i,j)$ 
and the Betti numbers of $\XX$ in case of 
$\mathbb Z_2$-homology.
The numbers $\mu_r(i,j)$ determine the 
numbers $\beta_r(i,j)$ for this filtration. 
Let $\XX=\{\sigma_1,\ldots,\sigma_n\}$ where ordering of the cells is 
topologically consistent and filtration compatible. 
An entry $M[i,j]$ in the incidence matrix
$M$ is set to the incidence number of $\sigma_i$ and $\sigma_j$, that is,
$M[i,j]:=I(\sigma_i,\sigma_j)$ where
$I(\sigma_i,\sigma_j)$ is $1$ if $\sigma_i$ is a face of codimension
$1$ of $\sigma_j$ and $0$ otherwise. 

Let $\mathcal U$ be the class of $n\times n$ upper triangular matrices 
with zero on diagonal and all entries zero or $1.$ 
The matrix $M$ is in this class. 
Given a matrix $M\in {\mathcal U}$, for any column  $j$, 
denote by $low(j)$ the largest $i$ so that $M[i,j]\ne 0.$  
If  the column $j$ has all entries zero $low(j)$ is not defined. 

The persistence algorithm uses column additions to transform  
the matrix $M$ within the class $\mathcal U$  with the additional property:
no two columns $j$ and $j'$ have $low(j)=low(j')$ if they are defined. 
A matrix with this property is called in ``reduced form". 
The algorithm works  by adding columns from left to right. 
Finally, the algorithm ends up with a matrix in the reduced  form. 
The matrix in the reduced form obtained from $M$ carries  explicit 
homological information about $\XX$ and the filtration.   
For example, the dimension of $H_r(\XX;\mathbb Z_2)$ is the 
cardinality of the set of zero columns $j$ in the reduced form
where $\sigma_j$ has dimension $r$
minus the number of nonzero columns $j'$ in the reduced form
where $\sigma_{j'}$ has dimension $r+1.$
To describe the number $\mu_r(i,j)$ 
we consider all pairs $(\sigma_k,\sigma_{k'})$  with  $low(k')=k$,
and observe that $\mu_r(i,j)$ is the number of such pairs 
which in addition satisfy 
\begin{enumerate} 
\item $k$  is an index corresponding to a cell in $\XX_i\setminus \XX_{i-1}$,
\item $k'$ is an index corresponding to a cell in $\XX_j\setminus \XX_{j-1}$,
\item dimension of $\sigma_{k}$ is exactly $r$. 
\end{enumerate}
\vskip .1in

Now we describe how one may use the persistence algorithm
to compute the numbers $\omega_r(\cdot,\cdot)$ that we defined for
simultaneous persistence.
Suppose $\XX$ is a cell complex and $\XX^\pm\subset \XX$ are
two subcomplexes with $A=\XX^+\cap \XX^-.$ 
Suppose $\XX^+$ and $\XX^-$ are equipped with finite filtrations 
$\XX^+_0\subseteq \XX^+ _1\subseteq \XX^+_2 \cdots \XX^+_{N^+}= \XX^+$ and
$\XX^-_0\subseteq \XX^- _1\subseteq \XX^-_2 \cdots \XX^-_{N^-}= \XX^-$
where $\XX^\pm_0=A.$

One can put the cells of $\XX$ in three groups I, II, III.   
The group I contains the cells of A, the group II the cells 
of $\XX^-\setminus A$ and the group III the cells of $\XX^+\setminus A.$ 
Suppose the cells are ordered according to condition A and condition B 
w.r.t. the filtration $ \XX_0\subseteq \XX_1\subseteq 
\XX_2 \cdots \XX_{N^-+N^+}=\XX$ defined by  
$\XX_i = \XX^-_i$ for $0\leq i\leq N^-$ and 
$\XX_i = \XX^-\cup \XX^+_{i-N^-}$ for   $N^-+1\leq i \leq N^-+ N^+.$
The incidence matrix $M$ is given by 
$
M=
\begin{Vmatrix}
A & B^- & C^+\\
0 & B & 0 \\
0 & 0 & C
\end{Vmatrix}
$
with $A$ being the incidence matrix for the complex $A$, the matrix 
$M^-= 
\begin{Vmatrix}
A & B^-\\
0 & B 
\end{Vmatrix}
,$
being the incidence matrix for $\XX^-,$ and  
$
M^+=
\begin{Vmatrix}
A & C^+\\
0 & C
\end{Vmatrix}
,$
being the incidence matrix for $\XX^+.$

The matrix $M$ is said to be  in {\it  relative reduced form} if $M^-$ 
and $M^+$ are in reduced form. By  running the persistence
algorithm first on $A$ and then continuing on $M^-$ and  $M^+$,
the relative reduced form can be achieved.

To compute the number $\omega_r(i,j)$, 
we consider all triples  
$(\sigma_k,\sigma_{k'}, \sigma_{k''})$  with  $low(k')= low(k'')=k$,
and observe that $\omega_r( i, j)$ is the number of such triples 
which satisfy in addition 
\begin{enumerate} 
\item $k$ is an index corresponding to a $r$-cell in $A$
\item $k'$  is an index corresponding to a cell in $\XX_i\setminus \XX_{i-1}$,
\item $k''$ is an index corresponding to a cell in $\XX_{j +N^-}\setminus \XX_{j+N^- -1}$.
\end{enumerate}

Next, we describe how we adapt the standard persistence algorithm to
compute the relevant numbers for level persistence. First, we observe the
following. If $\XX$ is a simplicial complex, $f:\XX\to \mathbb R$ is 
a generic linear map,  
and $t_0< t_1 ,\cdots t_i< \cdots $ are 
the values of $f$ on vertices of $\XX,$ then the filtration of $\XX$ by 
topological spaces $f^{-1}((-\infty, t_i])$ is not a filtration by 
simplicial sub complexes.  Apparently  the above algorithm can  
not be applied.  However it is possible to show that this topological 
filtration is homotopically the same (hence has the same persistence) 
as the filtration provided by the sub complexes $\XX_i$ consisting of 
the union of simplices $\sigma$ so that $f(\sigma)\subset (-\infty, t_i]$.


\paragraph{Computing level persistence for $\mathbb R$-valued maps.}
In the case of level persistence, the calculation of the 
numbers $l^\pm (i,j)$ and $\nu^\pm_r(i,j)$ for a tame map is reduced 
to the calculation of the numbers $\mu_r(\cdot, \cdot)$ 
for subspaces of $\XX,$ precisely  $\lx{t}$,  $\lx{t,\infty}$,
$\lx{-\infty, t}$equipped with the appropriate function , 
and $\lx{t_1, t_2}$ as indicated at the end of section~\ref{lpers-sec}.
In what follows we describe how to derive the matrix $M$ 
for $\lx{t}$, $\lx{t,\infty}$, $\lx{-\infty, t}$,  $\lx{t_1, t_2}$ 
from the  matrix $M$ for $\XX$.



Let $\XX$ be a simplicial complex and $f:\XX\to \mathbb R$ 
be a generic linear map.  
Observe that $\lx{t},$  $\rx{t}{\infty},$ $\rx{-\infty}{t},$ $\rx{t_1}{t_2}$ 
are cell complexes as indicated in the previous section and  
not simplicial complexes. However, the incidence structure
of their cells can be described  
in terms of the incidences of simplices in $\XX$.
This key observation allows us to construct the incidence
matrix of the cell complexes 
$\lx{t},$  $\rx{t}{\infty},$ $\rx{-\infty}{t},$ $\rx{t_1}{t_2}$ 
from the matrix $M(\XX)$ as detailed below.
For any simplex $\sigma$, let 
$[\sigma]$ denote the closed interval which is
the convex hull of the numbers 
$\{f(v), v \mbox{ a vertex of $\sigma$}\}$. 


\vskip .1in
\noindent
{\it Cell complex $\lx{t}.$}
For any $t$ introduce the matrix $\hat M(t)$ as the minor of $M(\XX)$ 
consisting of the rows and columns $i$ with $t\in \rm {int}\, 
[\sigma_i]$ and regard the simplex $\sigma_i$ as the 
cell $\hat \sigma_i$ with 
${\rm dim}\,(\hat \sigma_i)= {\rm dim}\,(\sigma_i)-1.$ 
If there are $\ell<n$ such simplices,
$\hat M(t)$ is an $\ell\times \ell$ matrix. 
Clearly, $\hat M(t)$ respects the order of the cells
induced by the order of simplices used for $M(\XX)$.  

\begin{itemize}
\item  If no vertex takes the value $t,$ 
then $M(\lx{t})= \hat M(t)$ is the incidence matrix of the 
cell complex $\lx{t}.$

\item  If there exists a vertex 
$v$ (unique since $f$ is generic) so that $f(v)=t$,
then the incidence matrix $M(\lx{t})$
is an $(\ell+1)\times (\ell+1)$ matrix with one additional row and column  
corresponding to the vertex $v$ viewed as a cell of dimension $0$ which 
should be indexed before all other cells. The entry for the
pair $(v,\hat\sigma_j)$ in $M(\lx{t})= 1$ if $\sigma_j$ is a $2-$simplex 
which has $v$ as a vertex and $0$ otherwise. 
\end{itemize}

\vskip .1in
\noindent
{\it Cell complexes $\lx{t, \infty},  \ \lx{- \infty, t}, \lx{t_1,t_2}.$}
To describe the incidence matrices which correspond to the cell 
complexes $\lx {t,\infty}$,
$\lx{-\infty, t}$, and $\lx{t_1,t_2}$ 
with $t_1< t_2$, consider the collection of simplices:
\begin{eqnarray*}
\mathcal X(t, \infty):= \{\sigma , [\sigma]\cap (t,\infty)\ne \emptyset\}\\
\mathcal X(-\infty, t):= \{\sigma , [\sigma]\cap (-\infty,t)\ne \emptyset\}\\
\mathcal X(t_1, t_2):= \{\sigma , [\sigma]\cap (t_1 , t_2)\ne \emptyset\}.
\end{eqnarray*}

Let the matrices $\hat M(t, \infty)$,  $\hat M(-\infty , t)$,
and $\hat M(t_1,t_2)$  be the minors of $M(\XX)$ whose 
rows /columns are the simplices in the sets 
$\mathcal X (t,\infty)$, $\mathcal X(-\infty,t)$,
and $\mathcal X(t_1,t_2)$ respectively. We keep the dimension 
of the cell same as that in the matrix $M(\XX)$.
The matrices $M(\lx{t,\infty})$ and $M(\lx{-\infty, t})$
are matrices of the form 

$$
\begin{Vmatrix}
\hat M(t) & B^+  \\
0 & \hat M(t, \infty)
\end{Vmatrix}
\mbox{ and }
\begin{Vmatrix}
\hat M(t)& B^-\\
0 & \hat M(-\infty , t)
\end{Vmatrix}
\mbox{ respectively},
$$
The matrix $B^+$ is a $s\times m$ matrix
where $\XX_t$ has $s$ cells and $\mathcal X(t,\infty)$
has $m$ simplices. We set $B^+[i,j]=1$ if the $i$-th cell in $\XX_t$
and the $j$-th simplex in $\mathcal X(t,\infty)$
are either $\hat \sigma$ and $\sigma$ respectively, or a vertex
of $\XX$ and a $1$-simplex incident to it respectively. All other
entries in $B^+$ are $0$. Define $B^-$ analogously 
by considering simplices in $\mathcal X(-\infty,t)$. 
The incidences among the cells in $\XX_{t,\infty}$ are of three
types, the ones among the cells in $\XX_t$, the ones among
the cells in $\XX_{t,\infty}$ that are not in $\XX_t$, and the
ones among the cells one of which is in $\XX_t$ and the other
is in $\XX_{t,\infty}$ but not in $\XX_t$. The matrices $\hat M(t) ,$ $ B^+ $
and $\hat M(t, \infty)$ capture these three types
of incidences respectively. Clearly, similar observations
can be made about the incidences among the cells in
$\XX_{-\infty,t}$.
\vskip .1in
Let $\hat M(t_1,t_2)$ denote the minor of $M(\XX)$ whose
rows and columns correspond to the simplices in $\mathcal X(t_1,t_2)$.
The matrix $M(\lx{t_1,t_2})$ is
$$
\begin{Vmatrix}
M(\lx{t_1}) & 0 & B^+\\
0 & M(\lx{t_2}) &B^- \\
0 & 0 & \hat M(t_1,t_2)
\end{Vmatrix}
$$
where 
the matrix $B^+$ has the entries 
as defined in the previous case with $\XX_{t_1}$ and $\mathcal X(t_1,t_2)$
playing the roles of $\XX_t$ and $\mathcal X(t,\infty)$ respectively.
The matrix $B^-$ is analogous.

The order of the cells as suggested by the matrices 
$M(X_t)$, $M(\lx{-\infty,t}), M(\lx{t,\infty})$ and $M(\lx{t_1,t_2})$  
satisfies  Condition A
but not necessarily Condition B w.r.t. the filtration induced from 
the map $g$ (end of section~\ref{lpers-sec}). 
If necessary it  can be canonically changed as 
indicated at the beginning of this section.

\paragraph{Computing level persistence for a circled value map.}
Let $f:\XX\to \mathbb S^1$ be a generic linear map.  For simplicity in writing we suppose that  $f(\sigma)$ is an angular interval smaller than $\pi.$
We need the incidence structure of the simplices in
the covering space $\tilde \XX$ and the linear $\mathbb R$-valued 
map $\tilde f$. However, the space $\tilde \XX$ is infinite.
So, we compute only a finite subspace of $\tilde \XX$ which
is sufficient for computing relevant persistence numbers and the restriction of $\tilde f$  to these subspaces. 

We show how to construct a simplicial complex from $\XX$ that contains 
$\tilde \XX_{t, t+ 2\pi k}$. The constructions of
$\tilde \XX_{t-2\pi k, t}, \ \tilde \XX_{t-2\pi k_1, t+ 2\pi k_2}$
are analogous.

\begin{figure}[h]
\centerline{\input{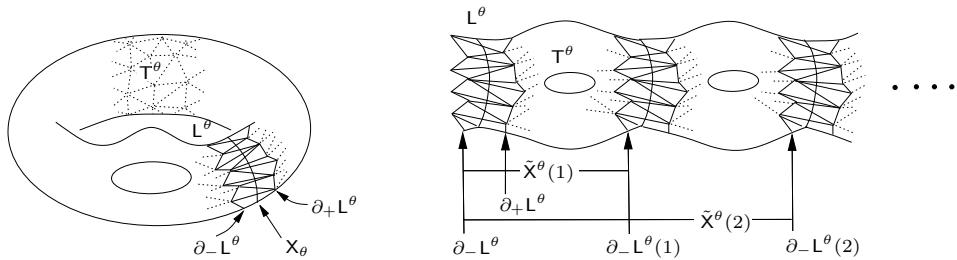}} 
\caption{Complex $\XX$ with level $\XX_\theta$ on left. Complexes
$\tilde\XX^\theta(1)$, $\tilde\XX^\theta(2)$, $\cdots$
providing a filtration of the total space $\tilde\XX$ on right.}
\label{fg:cover}
\end{figure}

Let $\mathcal X$ denote the set of simplices in $\XX$.
For any $\theta\in \mathbb S^1$, decompose $\mathcal X$ as a disjoint union 
$\mathcal X= \mathcal T^\theta \sqcup\mathcal L^\theta\sqcup \partial_-\mathcal L^\theta\sqcup \partial_+\mathcal L^\theta$
where (see Figure~\ref{fg:cover})
\begin{itemize}
\item $\mathcal L^\theta$ consists of the set of all 
simplices whose closure do intersect the level $\XX_\theta.$ 
Let $\LL^\theta$ be the simplicial complex generated by  
simplices in $\mathcal L^\theta,$ 

\item $\mathcal T^\theta$ is the set of simplices which do 
not belong to $\LL^\theta$.
Let $\TT^\theta$ denote the simplicial complex generated by the the 
simplices in $\mathcal T^\theta$ and consider $\TT^\theta \cap \LL^\theta$. 
This simplicial complex is the disjoint union of two simplicial
complexes $\partial_-\LL^\theta$ and $\partial_+\LL^\theta$ characterized by 
$f(\sigma)<\theta$ for $\sigma\in\partial _-\LL^\theta$ 
and $f(\sigma)>\theta$ for $\sigma\in\partial _+\LL^\theta$.

\item $\partial_\pm \mathcal L^\theta$  represent the 
simplices in $\partial_\pm \LL^\theta$.
\end{itemize}

Our purpose is to
build a collection of simplicial complexes which are equivalent with 
the space $\tilde \XX_{t, t+2\pi k}$
where $p(t)= \theta.$ 

Introduce a nested sequence of simplices 
$\tilde {\mathcal X}^\theta (0) \subseteq  
\tilde {\mathcal X}^\theta(1)\subseteq  \cdots 
\tilde {\mathcal X}^\theta(k)$ as follows.
Since we will repeat copies of each of the sets
$ \mathcal T^\theta$, $\mathcal L^\theta$, $\partial_-\mathcal L^\theta$,  
and $\partial_+\mathcal L^\theta$, let 
$ \mathcal T^\theta(n)$, $\mathcal L^\theta(n)$, 
$\partial_-\mathcal L^\theta(n)$,  and $\partial_+\mathcal L^\theta(n)$ 
denote their $n-$th copies respectively.
Taking ${\mathcal L}^\theta(0)={\mathcal L}^\theta$,
${\mathcal T}^\theta(0)={\mathcal T}^\theta$,
$\partial_+{\mathcal L}^\theta(0)=\partial_+{\mathcal L}^\theta$,
define inductively,
\begin{eqnarray*}
\tilde {\mathcal X}^\theta(0)&=& \partial_-{\mathcal L}^\theta\\
\tilde {\mathcal X}^\theta(n+1)&=& \tilde{\mathcal X}^\theta(n) 
\sqcup  \mathcal L^\theta(n) \sqcup\partial_+ \mathcal L^\theta (n) \sqcup \mathcal T^\theta(n) \sqcup \partial_- \mathcal L^\theta(n+1) 
\end{eqnarray*}

Taking $I_0(\sigma,\tau)=I(\sigma,\tau)$ and assuming $\tau$ to be a face
of $\sigma$ of codimension $1$,
the incidences among the simplices are described by 
\begin{eqnarray*}
I_{n+1}(\sigma,\tau)&=& I_{n}(\sigma, \tau) \mbox{ if 
$\sigma \in  \tilde{\mathcal X}^\theta(n)$ }
\\
I_{n+1}(\sigma,\tau)&=& I(\sigma, \tau) \mbox{ if $\sigma\in   \mathcal L^\theta(n)\sqcup \partial _+\mathcal L^\theta(n) \rm{\ viewed \ as \ simplices \ of} ~ \LL^\theta$ }
\\
I_{n+1}(\sigma,\tau)&=& I(\sigma, \tau) \mbox{ if $\sigma\in  \mathcal T^\theta(n)\sqcup \partial _-\mathcal L^\theta(n+1) \rm{\ viewed\  as \ simplices \ of}  ~ \TT^\theta$.}
\end{eqnarray*}
In all other cases $I_{n+1}(\sigma,\tau)=0$.
Notice that each $\tilde{\mathcal X}^\theta(i)$ forms a simplicial
complex $\tilde \XX^\theta(i)$ (Figure~\ref{fg:cover}).
To describe $\tilde f$ it suffices to provide its values on vertices.  We write  
$$\mathcal P = \partial _-\mathcal L^\theta\sqcup \mathcal L^\theta \sqcup  \partial _+\mathcal L^\theta \sqcup 
\mathcal T^\theta$$  
and $\mathcal P_0$ for the subset of vertices in $\mathcal P.$ We write $\mathcal P_0(n)$ for the $n-$th copy of $\mathcal P_0$
and define $\tilde f(n): \mathcal P_0(n) \to \mathbb R$ by $\tilde f(n):= \tilde f +2\pi n$ where 
$\tilde f= p^{-1}\cdot f$ with $p:[t, t+2 \pi)\to \mathbb S^1$ which
\footnote{$p$ is bijective and continuous, but $p^{-1}$ is not continuous} 
sends $t$ to $\theta$.
Once defined on vertices,
$\tilde f$ is extended by linearity to each simplex of the 
simplicial complex $\tilde \XX^\theta (n).$

Note  that  an order of the simplices of 
$\XX$ satisfying condition A induces an order on the simplices  
of $\mathcal T^\theta(n)$,
$ \partial_\pm \mathcal L^\theta(n)$ and $ \mathcal L^\theta(n)$ 
and by  juxtaposition an order on the simplices of 
$\tilde {\XX}^\theta(n)$ which continue to satisfy condition A. 
It implies that one can build a matrix
$M(\tilde\XX^\theta(n))$ which satisfies condition A
by juxtaposing the minors of $M(\XX)$ that represent $\mathcal L^\theta$,
$\mathcal \partial_{\pm}L^\theta$, $\mathcal T^\theta$,
and their copies in an appropriate order.
Note also that 
$\tilde {\XX}^\theta(n)$ is a sub complex of $\tilde \XX$ and therefore 
the restriction of $\tilde f$ provides tame maps on each of these spaces. 
The columns and rows of $M(\tilde\XX^\theta(n))$ can be reordered so that
they become filtration compatible with $\tilde f$ or $g$
by the method  indicated at the beginning of this section.

\cancel{
\Dan { I do not think that the remaining part is necessary}

\paragraph{Modifying the filtration.}  
The data for a circle valued linear
map $f:\XX\to \mathbb S^1$ is presented with an incidence matrix $M=M(\XX)$.
Assuming a total ordering among the simplices $\{\sigma_1,\cdots,\sigma_n\}$
one sets $M[i,j]=1$ if $\sigma_i$ is a face of codimension one
of $\sigma_j$ and $M[i,j]=0$ otherwise.
We also have $f$-values at the vertices of $\XX$
from which $f(\sigma)$ for any $\sigma\in \XX$ can be computed.

Given such an input, 
we compute an $n\times n$ matrix $M_0$ for the 
set $\tilde{\mathcal X}^{\theta}(0)$
from $M$ as follows: 
\begin{eqnarray*}
M_0[i,j]= \left\{\begin{array}{ll}
0 & \mbox{ if $\sigma_i\in \partial_+ \mathcal L^\theta$
 and $\sigma_j\in \mathcal L^\theta$}\\
M[i,j] & \mbox{otherwise}
\end{array}
\right.
\end{eqnarray*} 
Next, we compute the matrix $M_k$ for $\tilde{\mathcal X}^\theta(k)$
inductively:
$$
M_k=
\begin{Vmatrix}
M_{k-1}& D\\
D^T & M_0
\end{Vmatrix}
\mbox{ where $nk\times n$ matrix}
D[i,j]=\left\{
\begin{array}{ll}
1 & \mbox{ if $\sigma_i\in \partial_+\mathcal L^\theta(k)$
 and $\sigma_j\in \mathcal L^\theta(k-1)$}\\
0  & \mbox{otherwise}
\end{array}
\right.
$$
The only task left is that $M_k$ may not have the simplices
ordered according to an admissible filtration respecting
$\tilde f:\tilde \XX\to \mathbb R$. For this
move the columns and rows of $M_k$ as follows.
First organize $M_k$ by a total order on the simplices
according to dimensions, that is, if $\sigma_i$ is a face
of codimension 1 of a simplex $\sigma_j$, then $i< j$.
Next, move a column corresponding to a simplex $\sigma$ to the right 
of the simplex $\sigma'$ where $\sigma'$ is the
rightmost simplex with $f(\sigma') < f(\sigma)$.
Of course, in the above comparison, it is understood
that $f$-values are taken on the original copies
of $\sigma$ and $\sigma'$ in $\XX$.
Similarly, adjust the rows by moving them down.
}

\paragraph{Algorithm for almost integral 1-cocycles.}
If $(\bold f, \alpha)$ is an almost integral 1-cocycle, 
hence given by $M(\XX)$ and  a real number $\bold f  (e)$ for each edge, 
we have to provide the map $f:\XX \to  \mathbb S^1$  hence a matrix 
$M(\XX)$ and an angle valued map on vertices and calculate the level 
persistence. For this purpose we choose a base point vertex $x$ 
and assign to it the angle value $0$. For any other vertex $y$ choose 
a sequence of consecutive vertices  $\{x= y_0, y_1, \cdots y_{L-1}, y_L=y\}$ 
and assign to $y$ the angle 
$(\sum_{0\leq i\leq L-1} \bold f(y_i, y_{i+1})) \  
mod \ \alpha) 2\pi/ \alpha.$  
The result is independent of the choice of $x$. We continue then as 
described above. The total number of copies of $\LL$ and $\TT$  
involved in the calculation of all persistence numbers 
described here is less than the 
cardinality of $\mathcal X.$

\newpage
\section*{Appendix}

In this Appendix we show that, for a tame map $f:\XX\to \mathbb R$,
the relevant level persistence  numbers determine the persistence 
numbers for both $f$ and $-f.$  We suppose $\XX$ is compact 
so $f$ is bounded from above and below. 

Let $\{t_1, t_2, \cdots  t_k\}$ be the critical values of $f$ and let $\{ t'_1 ,\cdots t_{k+1}'\}$ regular values so that 
$t'_1 < t_1 <t'_2 < t_2 <\cdots t'_k < \cdots t_k < t'_{k+1}.$ 
Consider $s_{2i}= t_i$ and $s_{2i-1}= t'_i.$   
Clearly  $\lx{s_1}=\emptyset$   and $\lx{s_{2k+1}}= \XX,$  
therefore $\lx{s_1, s_i}= \lx{ -\infty, s_i}$  
and $\lx{s_i, s_{2k+1}}= \lx{s_i, \infty}.$ Note that the tameness 
of $f$ implies that the inclusions $\XX_{s_i}\subset \XX_{s_i, s_{i+1}}$ 
for $i$ even
and $\XX_{s_{i+1}} \subset \XX_{s_i, s_{i+1}}$ for $i$ odd 
induce isomorphism in homology. 

The relevant level persistence numbers are $ l_r(s_i)$, 
$ l^\pm_r (s_i; \pm(s_{i\pm k}- s_{i})$ and 
$e_r(s_i ; s_i-s_{i-k'}, s_{i+k''}- s_{i})$
while the relevant persistence numbers for $f$ are  
$\kappa_r(i)= \dim H_r(\lx{-\infty,s_i})$  and 
\newline $\kappa_r(s_i, s_{i+k})= \dim \ker (H_r (\XX_{-\infty, s_i}) \to H_r(\XX_{-\infty, s_{i+k}}) )$ while for $-f$ are  $\dim H_r(\lx{-s_i, \infty})$  and 
\newline  $\dim  \ker (H_r(\XX_{ -s_i,\infty}) \to H_r(\XX_{ -s_{i+k}, \infty})).$ 

To show that the numbers $l_r, l^\pm_r , e_r$  determine the  numbers $\kappa _r$ we will use  the algebraic  
proposition below:

\begin{proposition}\label {P100}

Let $$0\to A_N\overset {\alpha_N}\to B_N\overset{\beta_N}\to \cdots \overset {\delta_{n+1}}\to A_n\overset {\alpha_n}\to B_n \overset {\beta_n}\to C_n\overset {\delta_n}\to A_{n-1}\overset {\alpha_{n-1}} \to \cdots \to A_0\overset {\alpha_0}\to B_0\overset {\beta_0}\to C_0\overset {\delta_0=0}\to 0 \ \  \ (\ast) $$
$0 \leq n \leq N$ be an exact sequence of vector spaces. Then any three independent  
collections of numbers $\{\dim A_n\}, \{\dim B_n\},$
$ \{\dim C_n\}, \{\dim (\ker \alpha_n)\}, \{\dim (\ker \beta_n)\},
\{ \dim (\ker \delta_n)\} $ determine the other three.
(Here independent means  not constrained by the 
obvious equalities (1), (2), (3)  below.
\end{proposition}

We have  the following two long exact sequences:

\begin{enumerate}
\item 
The Mayer- Vietoris sequence associated to 
$\lx{t,t''}= \lx{t,t'}\cup \lx{t', t''}, \ \lx{t,t'}\cap \lx{t', t''}= \lx{t'},$ $t\leq t'\leq t".$
$$ \cdots \overset {\alpha_n}\to H_n(\lx{t,t'} )\oplus H_n(\lx{t',t''} ) \overset {\beta_n}\to H_n(\lx {t,t''}) \overset {\delta_n}\to 
 H_{n-1} (\lx{t'}) \overset {\alpha_{n-1}}\to H_{n-1}(\lx{t,t'} )\oplus H_{n-1}(\lx{t',t''} ) \overset {\beta _{n-1}}\to \cdots.$$

\item The long exact sequences  of the pairs $\lx{t, t'}\subset \lx{t, t''}$ and  $\lx{t', t''}\subset \lx{t, t''}, \ t\leq t' \leq t''.$
$$ (I)  \hskip .2in \cdots \overset {\delta_{n+1}}\to H_n(\lx{t,t'} ) \overset {\iota_n}\to H_n(\lx {t,t''}) \overset {j_n}\to 
 H_{n} (\lx{t,t''}, \lx{t,t'} ) \overset {\delta_{n}}\to H_{n-1}(\lx{t,t'} ) \overset {\iota _{n-1}}\to \cdots.$$

$$ (II)\hskip .2in \cdots \overset {\delta_{n+1}}\to H_n(\lx{t',t''} ) \overset {\iota_n}\to H_n(\lx {t,t''}) \overset {j_n}\to 
 H_{n} (\lx{t,t''}, \lx{t',t''} ) \overset {\delta_{n}}\to H_{n-1}(\lx{t',t''} ) \overset {\iota _{n-1}}\to \cdots.$$
\end{enumerate}

For $t\leq t'\leq t''$, we also have by excision theorem

$$ (III) \hskip .1in   H_{n} (\lx{t,t''}, \lx{t,t'} )=  H_{n} (\lx{t',t''}, \lx{t'} ),  \  H_{n} (\lx{t,t''}, \lx{t',t''} )=  H_{n} (\lx{t,t'}, \lx{t'} ).$$

\vspace{0.1in}
\noindent
\begin{proof}[Proposition \ref{P100}]:
Note first that the exactness $(\ast)$ implies  the following equalities 

\begin{equation}
\begin{aligned}
 (1)  \dim A_n= &\dim \ker \alpha_n +\dim \ker \beta_n\\
 (2) \dim B_n= &\dim \ker\beta_n + \dim \ker \delta_n\\
 (3) \dim C_n= &\dim \ker \delta_n + \dim \ker \alpha_{n-1}
\end{aligned}
\end{equation}
Then we obtain 
\begin{itemize}
\item
$\dim A_n, \dim B_n, \dim \ker \alpha_n,$  $0\leq n \leq N$  in view of  (1), (2) ,(3) determine  $\dim \ker \beta_n, \dim \ker \delta_n,$  and $\dim C_n$ ,$0\leq n \leq N.$
\item
$\dim A_n, \dim B_n, \dim \ker \beta_n,$   $0\leq n \leq N$  in view of  (2), (1), (3) determine $\dim \ker \delta_n,$  $\dim \ker \alpha_n$ and $\dim C_n$ ,$0\leq n \leq N.$ 
\item
$\dim A_n, \dim C_n, \dim \ker \alpha_n,$  $0\leq n \leq N$  in view of  (1), (3), (2) determine $\dim\ker \beta_n, \dim \ker \delta_n,$  and $\dim B_n$ ,$0\leq n \leq N.$ 
\item
$\dim A_n, \dim C_n, \dim \ker \beta_n,$  $0\leq n \leq N$ in view of  (1), (3), (2) determine , $\dim\ker \alpha_n, \dim \ker \delta_n$  and $\dim B_n$ ,$0\leq n \leq N.$ 
\item
$\dim \ker \alpha_n, \dim \ker \beta_n, \dim B_n,$ $0\leq n \leq N$  in view of  (1), (2), (3)determine  $\dim A_n, \dim \ker \delta_n$  and $\dim C_n$ ,$0\leq n \leq N.$ 
\item
$\dim \ker \alpha_n, \dim \ker \delta_n, \dim B_n,$ $0\leq n \leq N$  in view of  (3), (3), (2) determine $\dim C_n, \dim \ker \beta_n$  and $\dim A_n$ ,$0\leq n \leq N.$ 
\item
$\dim \ker \alpha_n, \dim \ker \delta_n, \dim A_n,$ $0\leq n \leq N$   in view of  (1), (2) , (3)determine  $\ker \dim \beta_n, \dim B_n$  and $\dim C_n$ ,$0\leq n \leq N.$ 
\item
$\dim \ker \alpha_n, \dim \ker \beta_n\dim \ker \delta_n,$  $0\leq n \leq N$ in view of  (1), (2) , (3) determine  
$\dim A_n,$ $\dim B_n,$ and $\dim C_n,$  $0\leq n \leq N$ 
\item
$\dim A_n,$ $\dim \ker \alpha_n, \dim \ker \delta_n$  ,$0\leq n \leq N$ in view of  (1), (2) , (3) determine  
$\dim \ker \beta_n,$ $ \dim B_n$  and $\dim C_n ,$ $0\leq n \leq N.$ 
 \item 
$\dim A_n,$ $\dim B_n,$ $\dim C_n,$  $0\leq n \leq N$ in view of  (1), (2) , (3) and the fact that $\dim \ker \alpha_N= 0$ 
determine  
$\dim\ker\alpha_n, $\ $dim\ker \beta_n ,$ $\dim \ker \delta_n,$ $0\leq n \leq N.$
\end{itemize}
All other possible situations  can be recovered from these cases. 
\end{proof}

To conclude that the level persistence numbers determine the persistence numbers we proceed as follows.
Use the equality  $\dim H_r(\XX_{s_i, s_{i+1}})= \dim H_r(\lx{s_i})$ if $i$ even  and  $\dim H_r(\XX_{s_i, s_{i+1}})= \dim H_r(\lx{s_{i+1}})$ if $i$ odd
which follow from tameness,
and  apply  Proposition \ref{P100} to the exact sequence (1). 
In this way we  derive from 
the numbers $l_r(s_i)$ and $e_r(s_i; s_i-s_{i-k'}, s_{i+k''}-s_i)$  
the numbers  
$\dim H_r (\lx{s_i, s_{i+k}})$ for all $i, k.$
\vskip .1in 
Apply Proposition \ref {P100} to the long exact sequence 2(I)  for $t=t',$ 
(resp. 2(II)  for $t'=t''$), and in view of  (III)  derive from the 
numbers $l_r(s_i)$ and  $l^+_r(s_i; s_{i+k}-s_i)$ 
(resp. $l^-_r(s_i; s_i-s_{i-k})$)
the numbers \newline $\dim H_r(\lx{s_i, s_{i+r+k}}, \lx{s_{i+r}, s_{i+r+k}})$ 
(resp. $\dim H_r(\lx{s_i, s_{i+r+k}}, \lx{s_i, s_{i+r}})$).
\vskip .1in
Apply Proposition \ref {P100} to the long exact sequence 2(I) 
(resp. 2(II)) and derive from  $\dim H_r (\lx{s_i, s_{i+k}})$  
$\dim H_r(\lx{s_i, s_{i+r+k}}, \lx{s_i, s_{i+r}})$ (resp. 
$\dim H_r(\lx{s_i, s_{i+r+k}}, \lx{s_{i+r}, s_{i+r+k}})$)
the numbers   $\dim \ker (H_r(\lx{s_i, s_{i+r}})\to    
H_r(\lx{s_i, s_{i+r+k}}))$
(resp. \  $\dim \ker (H_r(\lx{s_{i+r}, s_{i+r+k}})\to   
H_r(\lx{s_i, s_{i+r+k}}))$). 
Taking $s_i= s_1$ (resp. $s_{i+k+r}=s_{2k+1}$)
one obtains the persistence numbers for $f$ and $-f.$

\newpage
\begin{thebibliography}{99}

\bibitem{BD11} D. Burghelea and T. K. Dey. Persistence for circle valued maps.
(arxiv:1104.5646v1), 2011.

\bibitem{CSD09} G. Carlsson and V. de Silva and D. Morozov. 
Zigzag persistent homology and real-valued functions.
{\em Proc. 25th Annu. Sympos. Comput. Geom.}
(2009), 247--256.

\bibitem{CZ09}
G. Carlsson and A. Zomorodian.
The Theory of Multidimensional Persistence. with Gunnar Carlsson.
{\em Discrete Comput. Geom.} {\bf 42} (2009), 71--93. 

\bibitem{CEH07} D. Cohen-Steiner, H. Edelsbrunner, and J. Harer.
Stability of persistence diagrams. {\em Discrete Comput. Geom.}
{\bf 37} (2007), 103-120.

\bibitem{CEH09} D. Cohen-Steiner and H. Edelsbrunner and J. Harer.
Extending persistence using Poincar\'{e} and Lefschetz duality.
{\em Found. Comput. Math.} {\bf 9 (1)} (2009), 79--103.

\bibitem{CEM06} D. Cohen-Steiner, H. Edelsbrunner, and D. Morozov.
Vines and vineyards by updating persistence in linear time.
{\em Proc. 22nd Annu. Sympos. Comput. Geom.} (2006), 119--134.

\bibitem{DW07} T. K. Dey and R. Wenger. Stability of critical points 
with interval persistence.
{\em Discrete Comput. Geom.} {\bf 38} (2007), 479--512. 

\bibitem{ELZ02} H. Edelsbrunner, D. Letscher, and A. Zomorodian.
Topological persistence and simplification. {\em Discrete
Comput. Geom.} {\bf 28} (2002), 511--533.

\bibitem{EH} H. Edelsbrunner  and  J. L. Harer.
Computational Topology, An Introduction {\em book, AMS} 

\bibitem{FL10} P. Frosini and C. Landi. Stability of multidimensional
persistent homology with respect to domain perturbations.
{\em arxiv: 1001.1078v2 [math.AT]}, May 2010.

\bibitem{Hatcher}
A. Hatcher. Algebraic Topology. Cambridge U. Press, New York, 2002.

\bibitem{Novikov} S.P. Novikov. 
Quasiperiodic structures in topology. {\em Topological methods
in modern mathematics, Proc. Sympos. in honor of John
Milnor's sixtieth birthday, SUNY, Stony
Brook, New York, 1991}, eds. L. R. Goldberg and A. V. Phillips, Publish or
Perish, Houston, TX, 1993.


\bibitem{SJ09} V. de Silva and M. Vejdemo-Johansson. 
Persistent cohomology and circular coordinates.
{\em Proc. 25th Annu. Sympos. Comput. Geom.} (2009), 227-236.

\bibitem{ZC05} A. Zomorodian and G. Carlsson. Computing persistent 
homology. {\em Discr. Comput. Geom.} {\bf 33} (2005), 249--274.


\bibitem{JLY} Xiaoye Jiang, Lek-Heng Lim, Yuan Yao and Yinyu Ye. 
Statistical Ranking and Combinatorial Hodge Theory. 
{\em arxiv.org/abs/0811.1067.} 

\bibitem{YY} Yuan Yao. Combinatorial Laplacians and Rank Aggregation.
{\em 6th International Congress of Industrial and Applied Mathematics (ICIAM), 
Mini Symposium: 
Novel Matrix Methods for Internet Data Mining. Zurich, Switzerland}, 2007.
\end{thebibliography}
\end{document}